\documentclass[11pt]{amsart}
 \usepackage{amssymb,amsfonts,amscd,amsbsy}
\usepackage[mathscr]{eucal}
\usepackage{url}

\newtheorem{thm}{Theorem}[section]

\newtheorem{prop}[thm]{Proposition}

\theoremstyle{definition}

\begin{document}

\title[A supercongruence for generalized Domb numbers]
{A supercongruence for generalized Domb numbers} 
 
\author{Robert Osburn and Brundaban Sahu}

\address{School of Mathematical Sciences, University College Dublin, Belfield, Dublin 4, Ireland}

\address{IH{\'E}S, Le Bois-Marie, 35, route de Chartres, F-91440 Bures-sur-Yvette, FRANCE}

\address{School of Mathematical Sciences, National Institute of Science Education and Research, Bhubaneswar 751005, India}

\email{robert.osburn@ucd.ie, osburn@ihes.fr}

\email{brundaban.sahu@niser.ac.in}


\date{\today}

\begin{abstract}
Using techniques due to Coster, we prove a supercongruence for a generalization of the Domb numbers. This extends a recent result of Chan, Cooper and Sica and confirms a conjectural supercongruence for numbers which are coefficients in one of Zagier's seven ``sporadic" solutions to second order Ap{\'e}ry-like differential equations.
\end{abstract}

\maketitle

\section{Introduction}

It is now well-known that the {\it Ap{\'e}ry numbers}

$$
A(n):=\sum_{k=0}^{n} \binom{n}{k}^2 \binom{n+k}{k}^2
$$

\noindent play a crucial role in the irrationality proof of $\zeta(3)$, satisfy many interesting congruences and are related to modular forms. For example, Gessel \cite{gessel} showed that

\begin{equation} \label{ges}
A(np) \equiv A(p) \pmod{p^3}
\end{equation}

\noindent for any prime $p>3$, while if 

$$
F(z)= \frac{\eta^{7}(2z) \eta^{7}(3z)}{\eta^{5}(z) \eta^{5}(6z)} \quad \text{and} \quad t(z)= \Biggl( \frac{\eta(6z) \eta(z)}{\eta(2z) \eta(3z)} \Biggr)^{12},
$$

\noindent then by a result of Peters and Stienstra \cite{ps}, we have

$$
F(z) = \sum_{n=0}^{\infty} A(n) \, t^{n}(z).
$$

\noindent Here $\eta(z)$ is the Dedekind eta-function. Since then, there have been several papers which study arithmetic properties of coefficients of power series expansions in $t$ of modular forms where $t$ is a modular function (see \cite{beukers1}, \cite{ccs}, \cite{ckko}, \cite{jv}, \cite{os}, \cite{ossuper}, \cite{v}, \cite{z}). 

Our interest is in the sequence of numbers given by  

\begin{equation*}
D(n):= \sum_{k=0}^{n} \binom{n}{k}^{2} \binom{2k}{k} \binom{2(n-k)}{n-k}.
\end{equation*}

\noindent The first few terms in the sequence of {\it Domb numbers} $\{ D(n) \}_{n \geq 0}$ are as follows:

\begin{center}
$1$, $4$, $28$, $256$, $2716$, $31504$, $\dotsc$
\end{center}

\noindent This ubiquitous sequence (see $A002895$ of Sloane \cite{sloane}) not only arises in the theory of third order Ap{\'e}ry-like differential equations \cite{avsz}, odd moments of Bessel functions in quantum field theory \cite{bbbg}, uniform random walks in the plane \cite{bnsw}, new series for $1/\pi$ \cite{ccl}, interacting systems on crystal lattices \cite{domb} and the enumeration of abelian squares of length $2n$ over an alphabet with $4$ letters \cite{rs}, but if 

$$
G(z)= \frac{\eta^{4}(z) \eta^{4}(3z)}{\eta^{2}(2z) \eta^{2}(6z)} \quad \text{and} \quad s(z)= \Biggl( \frac{\eta(2z) \eta(6z)}{\eta(z) \eta(3z)} \Biggr)^6,
$$

\noindent then (see \cite{ccl})

$$
G(z) = \sum_{n=0}^{\infty} (-1)^n D(n) \, s^{n}(z). 
$$

Motivated by (\ref{ges}), Chan, Cooper and Sica \cite{ccs} recently proved the congruence

\begin{equation} \label{ccsd}
D(np) \equiv D(p) \pmod{p^3}.
\end{equation}

The purpose of this short note is to prove a supercongruence for the {\it generalized Domb numbers}. Recall that the term {\it supercongruence} refers to congruences that are stronger than those suggested by formal group theory (for recent developments in this area, see \cite{long}, \cite{mccarthy}, \cite{zudilin}). For integers $A$, $B$ and $C \geq 1$, let

\begin{equation} \label{gendomb}
D(n, A, B, C) := \sum_{k=0}^{n} \binom{n}{k}^{A} \binom{2k}{k}^B \binom{2(n-k)}{n-k}^C.
\end{equation}

\noindent Our main result is the following.

\begin{thm} \label{main} Let $A$, $B$ and $C$ be integers $\geq 1$ and $p>3$ be a prime. For any integers $m$, $r \geq 1$, we have

\begin{equation*} 
D(mp^r, A, B, C) \equiv D(mp^{r-1}, A, B, C) \pmod{p^{3r}}
\end{equation*}

\noindent if $A \geq 2$. 
\end{thm}

Note that Theorem \ref{main} recovers (\ref{ccsd}) in the case $A=2$, $B=C=1$, $r=1$ and generalizes a numerical observation in Section 3 of \cite{os} (see case (xii) in Table 3). The method of proof for Theorem \ref{main} is due to Coster in his influential Ph.D. thesis \cite{cos}. Namely, one expresses the summands in (\ref{gendomb}) as products $g_{AB}(X,k)$ and $g_{AB}^{*}(X,k)$ (see Section 2), then utilizes the combinatorial features of these products. One then writes (\ref{gendomb}) as two sums, one for which $p \mid k$ and the other for which $p \nmid k$. In the case $p \nmid k$, the sum vanishes modulo an appropriate power of $p$ while for $p \mid k$, the sum reduces to the required result. This strategy not only leads to a generalization of (\ref{ges}) (see Theorem 4.3.1 in \cite{cos}), but can be used to prove supercongruences for other similar sequences \cite{ossuper}. Additionally, a proof similar to that of Theorem 1.1 can be employed to show

$$
D(mp^r, 1, 1, 1) \equiv D(mp^{r-1}, 1, 1, 1) \pmod{p^{2r}},
$$ 

\noindent thereby confirming another conjectural supercongruence in Section 3 of \cite{os} (see case (ix) in Table 2). The details are left to the interested reader. The numbers $D(n,1,1,1)$ are coefficients in one of Zagier's seven ``sporadic" solutions (see $\# 10$ in Table 1 of \cite{z} or the modular parameterization given by Case E in Table 3 of \cite{z}) to a general family of second order Ap{\'e}ry-like differential equations. Our hope is that the present note will inspire others to further explore the techniques in \cite{cos}. In Section 2, we recall the relevant properties of the products $g_{AB}(X,k)$ and $g_{AB}^{*}(X,k)$ and then prove Theorem \ref{main}.

\section{Proof of Theorem \ref{main}}

We first recall the definition of two products and one sum and list some of their main properties. For more details, see Chapter 4 of \cite{cos}. For integers $A$, $B \geq 0$, $k$, $j \geq 1$ and $X$ and for a fixed prime $p>3$, we define

\begin{equation*} 
\displaystyle g_{AB}(X, k)=\prod_{i=1}^{k} \Biggl( 1 - \frac{X}{i} \Biggr)^{A} \Biggl(1 + \frac{X}{i} \Biggr)^{B},
\end{equation*}

\begin{equation*}
\displaystyle  g_{AB}^{*}(X, k)=\prod_{\substack{i=1 \\ p \nmid i}}^{k} \Biggl( 1 - \frac{X}{i} \Biggr)^{A} \Biggl(1 + \frac{X}{i} \Biggr)^{B},
\end{equation*}
 
\noindent and 

\begin{equation*}
\displaystyle S_{j}(k)=\sum_{\substack{i=1 \\ p \nmid i}}^{k} \frac{1}{i^j}.
\end{equation*}

The following proposition (see Lemmas 4.2.1 and 4.2.5 in \cite{cos}) provides some of the main properties of $g_{AB}(X, k)$, $g_{AB}^{*}(X,k)$ and $S_{j}(k)$. 

\begin{prop} \label{prop}
For any integers $A$, $B \geq 0$, $X \in \mathbb{Z}$ and integers $m$, $k$, $r \geq 1$, we have
\begin{enumerate}

\item[(i)] $\displaystyle S_{j}(mp^r) \equiv 0 \pmod{p^r}$ for $j \not\equiv 0 \pmod {p-1}$,

\item[(ii)] $\displaystyle S_{2j-1}(mp^r) \equiv 0 \pmod{p^{2r}}$ for $j \not\equiv 0 \pmod{\frac{p-1}{2}}$,

\item[(iii)] $\displaystyle g_{AB}(pX, k) = g_{AB}^{*}(pX, k) g_{AB}(X, \bigl \lfloor \tfrac{k}{p} \bigr \rfloor)$,

\item[(iv)] $\displaystyle g_{AB}^{*}(X,k) \equiv 1 + (B-A) S_{1}(k)X + \tfrac{1}{2} \Bigl( (A-B)^2 S_{1}(k)^2 - (A+B) S_{2}(k) \Bigr) X^2  \pmod {X^3}$,

\item[(v)] $\displaystyle \binom{n}{k}^{A} \binom{n+k}{k}^{B} =(-1)^{Ak} \Bigl ( \frac{n}{n-k} \Bigr)^{A} g_{AB}(n,k)$.

\end{enumerate}
\end{prop}

We now prove Theorem \ref{main}. \\

\noindent {\it{Proof of Theorem \ref{main}.}}
We first note that it suffices to prove the result with $p \nmid n$, $p \nmid m$ where $m$, $n \geq 1$ are integers and $p>3$ is a prime. We now assume that $A \geq 2$ and $B \geq C \geq 1$. Recall that for integers $m$, $n$, $r \geq 1$ with $p \nmid n$, $p \nmid m$ and $s \geq 0$ with $s \leq r$, we have

\begin{equation} \label{ord}
ord_{p} \binom{mp^r}{np^s}^{A} = A(r-s).
\end{equation}

\noindent Also, by Lemma 2.2 in \cite{ossuper}, we have for a prime $p>3$ and integers $m \geq 0$, $r \geq 1$ 

\begin{equation} \label{reduce}
\binom{2mp^r}{mp^r} \equiv \binom{2mp^{r-1}}{mp^{r-1}} \pmod{p^{3r}}.
\end{equation}

\noindent Now, taking $j=2$ in (i), $j=1$ in (ii) and $X=mp^r$, $k=np^s$ in (iv) of Proposition \ref{prop}, we have 

\begin{equation} \label{key}
g_{AB}^{*}(mp^r, np^s) \equiv 1 \pmod{p^{r+2s}}
\end{equation} 

\noindent for any non-negative integers $m$, $n$, $r$ and $s$ with $s \leq r$. Letting $n=mp^r$, $k=np^s$, $A \to A+2C$, $B=0$ in (v) and $X=mp^{r-1}$, $k=np^s$ in (iii) of Proposition 2.1, we have, for $s \geq 1$,

\begin{equation} \label{step1}
\begin{aligned}
\binom{mp^r}{np^s}^{A+2C} & = (-1)^{(A+2C)np^s} \Biggl( \frac{mp^r}{mp^r - np^s} \Biggr)^{A+2C} g_{(A+2C)\,0}(mp^r,np^s) \\
& = (-1)^{Anp^{s-1}}  \Biggl( \frac{mp^{r-1}}{mp^{r-1} - np^{s-1}} \Biggr)^{A+2C} g_{(A+2C)\,0}^{*}(mp^r,np^s) g_{(A+2C)\,0}(mp^{r-1},np^{s-1}) \\
& = \binom{mp^{r-1}}{np^{s-1}}^{A+2C} g_{(A+2C)\,0}^{*}(mp^r, np^s).
\end{aligned}
\end{equation}

\noindent In the last step of (\ref{step1}), we have applied (v) of Proposition 2.1 with $n=mp^{r-1}$, $k=np^{s-1}$, $A \to A + 2C$ and $B=0$. Thus,

\begin{equation} \label{allterms}
\begin{aligned} 
\binom{mp^r}{np^s}^{A+2C} \binom{2np^s}{np^s}^{B-C} &  \binom{2mp^{r-1}}{2np^{s-1}}^{C} \\
& \\
&= \binom{mp^{r-1}}{np^{s-1}}^{A+2C} g_{(A+2C)\,0}^{*}(mp^r, np^s) \binom{2np^{s}}{np^s}^{B-C} \binom{2mp^{r-1}}{2np^{s-1}}^{C}.
\end{aligned}
\end{equation}

\noindent Similarly, letting $n=2mp^r$, $k=2np^s$, $A=C$, $B=0$ in (v) and $X=2mp^{r-1}$, $k=2np^s$ in (iii) of Proposition 2.1, we have

\begin{equation} \label{step2}
\begin{aligned}
\binom{2mp^r}{2np^s}^{C} & = (-1)^{2Cnp^s}  \Biggl( \frac{2mp^r}{2mp^r - 2np^s} \Biggr)^{C} g_{C0}(2mp^r, 2np^s) \\
& =  \Biggl( \frac{2mp^{r-1}}{2mp^{r-1} - 2np^{s-1}} \Biggr)^{C} g_{C0}^{*}(2mp^r, 2np^s) g_{C0}(2mp^{r-1}, 2np^{s-1}) \\
& =  \binom{2mp^{r-1}}{2np^{s-1}}^{C} g_{C0}^{*}(2mp^r, 2np^s).
\end{aligned}
\end{equation}

\noindent In the last step of (\ref{step2}), we have taken $n=2mp^{r-1}$, $k=2np^{s-1}$, $A=C$ and $B=0$ in (v) of Proposition 2.1. By (\ref{reduce}) and  (\ref{key}), we have 

\begin{equation} \label{g2}
\binom{2np^s}{np^s}^{B-C} \equiv \binom{2np^{s-1}}{np^{s-1}}^{B-C} \pmod{p^{3s}}
\end{equation}

\noindent and

\begin{equation} \label{g1}
g_{(A+2C)\,0}^{*}(mp^r, np^s) \equiv g_{C0}^{*}(2mp^r, 2np^s) \equiv 1 \pmod{p^{r+2s}}.
\end{equation}

\noindent For $r \geq s$, $A \geq 2$ and $C \geq 1$, we now claim that

\begin{equation} \label{pdivk}
\displaystyle \frac{\binom{mp^r}{np^s}^{A+2C} \binom{2np^s}{np^s}^{B-C}}{\binom{2mp^r}{2np^s}^{C}} \equiv \frac{\binom{mp^{r-1}}{np^{s-1}}^{A+2C} \binom{2np^{s-1}}{np^{s-1}}^{B-C}}{\binom{2mp^{r-1}}{2np^{s-1}}^{C}} \pmod{p^{3r}}.
\end{equation}

\noindent To see this, we first note that by (\ref{step2}) and (\ref{g1}), (\ref{g2}) and (\ref{g1}), we have

\begin{equation} \label{y1}
 \binom{2mp^{r-1}}{2np^{s-1}}^{C} = \binom{2mp^r}{2np^s}^{C} - \gamma p^{r+2s} \binom{2mp^{r-1}}{2np^{s-1}}^{C},
\end{equation}

\begin{equation} \label{y2}
\binom{2np^s}{np^s}^{B-C} = \binom{2np^{s-1}}{np^{s-1}}^{B-C} + \alpha p^{3s}
\end{equation}

\noindent and

\begin{equation} \label{y3}
g_{(A+2C)\,0}^{*}(mp^r, np^s) = 1 + \beta p^{r+2s}
\end{equation}

\noindent for some $\gamma$, $\alpha$ and $\beta \in \mathbb{Z}$. After substituting (\ref{y1})--(\ref{y3}) into the right hand side of (\ref{allterms}) and multiplying, we consider the following seven terms:

\begin{enumerate}

\item[(a)] $\displaystyle p^{r+2s} \binom{mp^{r-1}}{np^{s-1}}^{A+2C} \binom{2np^{s-1}}{np^{s-1}}^{B-C} \binom{2mp^{r-1}}{2np^{s-1}}^{C}$ \\

\item[(b)] $\displaystyle p^{3s} \binom{mp^{r-1}}{np^{s-1}}^{A+2C} \binom{2mp^r}{2np^s}^{C}$ \\

\item[(c)] $\displaystyle p^{r+5s} \binom{mp^{r-1}}{np^{s-1}}^{A+2C} \binom{2mp^{r-1}}{2np^{s-1}}^{C}$ \\

\item[(d)] $\displaystyle p^{r+2s} \binom{2mp^r}{2np^s}^{C} \binom{2np^{s-1}}{np^{s-1}}^{B-C} \binom{mp^{r-1}}{np^{s-1}}^{A+2C}$ \\

\item[(e)] $\displaystyle p^{2r+4s} \binom{2np^{s-1}}{np^{s-1}}^{B-C} \binom{mp^{r-1}}{np^{s-1}}^{A+2C}  \binom{2mp^{r-1}}{2np^{s-1}}^{C}$ \\

\item[(f)] $\displaystyle p^{r+5s} \binom{2mp^r}{2np^s}^{C} \binom{mp^{r-1}}{np^{s-1}}^{A+2C}$ \\

\item[(g)] $\displaystyle p^{2r + 7s}  \binom{2mp^{r-1}}{2np^{s-1}}^{C} \binom{mp^{r-1}}{np^{s-1}}^{A+2C}$. \\

\end{enumerate}

\noindent As $ord_{p}$ is at least $3r + C(r-s)$ in each of the cases (a)--(g) above and we have (\ref{ord}), (\ref{pdivk}) follows. Now, using the identity

$$
\displaystyle \binom{a-b}{c-d} \binom{b}{d} = \frac{\binom{a}{c} \binom{c}{d} \binom{a-c}{b-d}}{\binom{a}{b}},
$$

\noindent we have

$$
\displaystyle D(mp^r, A, B, C) = \binom{2mp^r}{mp^r}^{C} \,\, \sum_{k=0}^{mp^r} \frac{\binom{mp^r}{k}^{A+2C} \binom{2k}{k}^{B-C}}{\binom{2mp^r}{2k}^{C}}.
$$

\noindent We now split $D(mp^r, A, B, C)$ into two sums, namely

$$
D(mp^r, A, B, C) = \binom{2mp^r}{mp^r}^{C} \,\, \sum_{\substack{k=0 \\ p \nmid k}}^{mp^r} \frac{\binom{mp^r}{k}^{A+2C} \binom{2k}{k}^{B-C}}{\binom{2mp^r}{2k}^{C}} + 
\binom{2mp^r}{mp^r}^{C} \,\, \sum_{\substack{k=0 \\ p \mid k}}^{mp^r} \frac{\binom{mp^r}{k}^{A+2C} \binom{2k}{k}^{B-C}}{\binom{2mp^r}{2k}^{C}}. 
$$

\noindent Since $A \geq 2$, $B \geq C \geq 1$, the first sum vanishes modulo $p^{3r}$ using (\ref{ord}) and the result then follows from reindexing the second sum and applying (\ref{reduce}) and (\ref{pdivk}). A similar argument holds in the case $A \geq 2$, $C > B \geq 1$ upon noting that

$$
\displaystyle \frac{\binom{mp^r}{k}^{A+2B} \binom{2(mp^r - k)}{mp^r - k}^{C-B}}{\binom{2mp^r}{2k}^{B}} \equiv 0 \pmod{p^{3r}}
$$

\noindent if $p \nmid k$ and 

$$
\displaystyle \frac{\binom{mp^r}{np^s}^{A+2B} \binom{2(mp^r - np^s)}{mp^r - np^s}^{C-B}}{\binom{2mp^r}{2np^s}^{B}} \equiv \frac{\binom{mp^{r-1}}{np^{s-1}}^{A+2B} \binom{2(mp^{r-1} - np^{s-1})}{mp^{r-1} - np^{s-1}}^{C-B}}{\binom{2mp^{r-1}}{2np^{s-1}}^{B}} \pmod{p^{3r}}
$$

\noindent if $p \mid k$. \qed

\section*{Acknowledgements} The first author would like to thank the Institut des Hautes {\'E}tudes
Scientifiques for their support during the revision of this paper while the second author would like to thank the UCD School of Mathematical Sciences for their hospitality during the preparation of this paper. The authors were partially funded by Science Foundation Ireland 08/RFP/MTH1081. Finally, the authors thank the referee for their careful reading of the paper, many helpful comments and patience.

\end{document}